\documentclass[reqno, 11pt]{amsart}
\usepackage[utf8]{inputenc}
\usepackage{amsmath}
\usepackage{amsthm}
\usepackage{amssymb}
\usepackage{mathtools}
\usepackage{booktabs}
\usepackage{amsxtra}
\usepackage{amsbsy}
\usepackage{comment}
\usepackage{mathrsfs}
\usepackage{longtable}

\usepackage[english]{babel}
\usepackage{blindtext}
\usepackage{url}
\usepackage{bbm}
\usepackage{euscript}
\usepackage{pifont}
\usepackage{hyperref}
\usepackage{wasysym}
\hypersetup{
    colorlinks=true,
    linkcolor=blue,
    filecolor=magenta,
    urlcolor=cyan,
}
\usepackage{graphicx}
\usepackage{epstopdf}

    \usepackage[
    top    = 2.56cm,
    bottom = 2.56cm,
    left   = 2.61 cm,
    right  = 2.61  cm]{geometry}

\newtheorem{thm}{Theorem}[section]
\newtheorem{lem}[thm]{Lemma}

\newtheorem*{conjecture*}{Conjecture}
\newtheorem*{thm*}{Theorem}

\theoremstyle{remark}
\newtheorem*{remark}{Remark}

\theoremstyle{definition}

\newcommand{\N}{\mathbb{N}}
\newcommand{\Z}{\mathbb{Z}}

\newcommand{\Alt}{{\operatorname{Alt}}}
\newcommand{\RNum}[1]{\uppercase\expandafter{\romannumeral #1\relax}}

\newcommand{\cod}{\normalfont{\mbox{cod}}}
\newcommand{\cd}{\normalfont{\mbox{cd}}}
\newcommand{\GL}{\normalfont{\mbox{GL}}}
\newcommand{\SL}{\normalfont{\mbox{SL}}}
\newcommand{\SU}{\normalfont{\mbox{SU}}}
\newcommand{\PSL}{\normalfont{\mbox{PSL}}}
\newcommand{\PSU}{\normalfont{\mbox{PSU}}}
\newcommand{\Aut}{\mbox{Aut}}
\newcommand{\irr}{\mbox{Irr}}

\newskip\aline \newskip\halfaline
\title{Huppert's analogue conjecture for $\PSL(3,q)$ and $\PSU(3,q)$}

\author[Y. Liu]{Yang Liu}
\address{School of Mathematical Science, Tianjin Normal University, Tianjin 300387, P.\,R. China}
\email{yliu@tjnu.edu.cn}

\author[Y. Yang]{Yong Yang}
\address{Department of Mathematics, Texas State University, 601 University Drive, San Marcos, TX 78666}
\email{yang@txstate.edu}

\begin{document}

\maketitle

\begin{abstract}
    Let $G$ be a finite group and $\chi\in \irr(G)$. The codegree of $\chi$ is defined as
    $\cod(\chi)=\frac{|G:\ker(\chi)|}{\chi(1)}$ and $\cod(G)=\{\cod(\chi) \ |\ \chi\in \irr(G)\}$ is called the set of codegrees of $G$. In this paper, we show that the set of codegrees of $\PSL(3,q)$ and $\PSU(3,q)$ determines the group up to isomorphism.
\end{abstract}

\section{Introduction}

 Let $G$ be a finite group and $\irr(G)$ be the set of all irreducible complex characters of $G$.
The concept of codegrees was originally introduced by Chillag and Herzog in \cite{Chillag1} where the codegree of $\chi$ was defined as $\frac{|G|}{\chi(1)}$ for character $\chi\in \irr(G)$. However, the definition was modified to $\cod(\chi)=\frac{|G:\ker(\chi)|}{\chi(1)}$ later by Qian, Wang, and Wei in \cite{Qian} so that there is no different meaning for $\cod(\chi)$  when $\chi$ is considered as a character in some quotient group of $G$. Because the relationship between codegree and degree are very close, it seems that we can also characterize the structure of groups by codegree, and there might be a corresponding codegree version for the problem of character degree. During the past few years the study about the character codegrees
has been very active and a lot of results about codegrees have been obtained, such as the relationship between the codegrees and the element orders, codegrees of $p$-groups, and groups with few codegrees, see for example \cite{Aliz,dl,i,ly1,ly2,q21}.

Denote by cd$(G) = \{\chi(1) \ |\ \chi \in \irr(G)\}$. Bertram Huppert raised the following conjecture (which has been verified for sporadic simple groups, alternating groups, and simple groups of Lie type with low rank): \\

\textbf{Huppert's Conjecture:} Let $H$ be any finite non-abelian simple group and $G$ a finite group such that cd$(G)$ = cd$(H)$. Then, $G \cong H \times A$, where $A$ is abelian.\\

Denote by cod$(G) = \{\cod(\chi) \ |\ \chi \in \irr(G)\}$. Recently people began to study a  similar conjecture related to the codegrees of a group. \\

\textbf{Codegree version of Huppert's conjecture:} Let $H$ be any finite non-abelian simple group and $G$ a finite group such that $\cod(G)=\cod(H)$. Then $G\cong H$.\\

This conjecture has been considered and shown to hold for $\PSL(2,q)$ in \cite{BahriAkh}. In \cite{Ahan}, the conjecture was proven for ${}^2B_2(2^{2f+1})$, where $f \ge 1$, $\PSL(3,4)$, $A_7$, and $J_1$. The conjecture also holds in the cases where $H$ is $M_{11}, M_{12}, M_{22}, M_{23}$, or $\PSL(3,3)$ by \cite{gkl}. In this paper, we continue to study this conjecture and establish the following result.

 \begin{thm}

  If $H$ is isomorphic to simple groups $\PSL(3,q)$ or $\PSU(3,q)$  and $G$ a finite group such that $\cod(G)=\cod(H)$. Then $G\cong H$.

\end{thm}

\section{Preliminary}

Firstly we give a list of simple groups with few degrees (\cite[Theorem 1.1]{Aziz}).

\begin{lem} Let $S$ be a non-abelian finite simple group. Then $|\mbox{cd}(S)|> 10$ or one of the following holds.
\begin{enumerate}
    \item $|\cd(S)|=4$ and $S= \PSL(2, 2^f)$, $f\geq 2$;

    \item $|\cd(S)|=5$ and $S=\PSL(2, p^f)$, $p\neq 2$ and $p^f > 5$;

    \item $|\cd(S)|=6$ and $S= {}^2B_2(2^{2f+1}), f\geq 1$ or $S = \PSL(3,4)$;

    \item $|\cd(S)|=7$ and $S= \PSL(3,3), A_7, J_1, M_{11}$;

    \item $|\cd(S)|=8$ and $S=\PSL(3,q)$ where $4< q\not\equiv 1\ (\bmod\ 3)$ or $S = \PSU(3,q)$ where $4 < q \not\equiv -1\ (\bmod\ 3)$ or $S=G_2(2)'$;

    \item $|\cd(S)|=9$ and $S=\PSL(3,q)$ where $4<q\equiv 1\ (\bmod\ 3)$ or $S=\PSU(3,q)$ where $4<q \equiv -1\ (\bmod\ 3)$.
\end{enumerate}
\end{lem}

In Table 1 we give the complete codegree sets for simple groups with few codegrees which will be used many times in the proof of the main results. One may check \cite{Conway, sf, su} for the details.

\begin{table}
\caption{Codegree set for some simple groups}

\begin{tabular}{ll}
 \hline
 Group $S$  & Codegree set $\cod(S)$ \\ [1pt] \hline

  $\PSL(2,k)$ ($k=2^f\geq 4$) & $\{1, k(k-1), k(k+1), k^2-1\}$\\
  \hline
$\PSL(2,k)$ ($k>5$ odd)  & $\{ 1, \frac{k(k-1)}{2}, \frac{k(k+1)}{2}, \frac{k^2-1}{2}, k(k-\epsilon(k))\}$, $\epsilon(k)=(-1)^{(k-1)/2}$ \\ [1pt] \hline

${}^2B_2(q)$ &  $\{1,\  (q-1)(q^2+1),\  q^2(q-1),\  2^{3f+2}(q^2+1),\ q^2(q-2r+1),$ \\ [1pt]
 $q=2r^2=2^{2f+1}$    &  $q^2(q+2r+1) \}$ \\ [1pt] \hline

 $\PSL(3,4)$ & $\{ 1,\  2^4{\cdot}3^2{\cdot}7,\ 2^6{\cdot}3^2,\ 2^6{\cdot}5,\ 2^6{\cdot}7,\ 3^2{\cdot}5{\cdot}7 \}$\\ [1pt] \hline

$\Alt_7$ & $\{1,\ 2^2{\cdot} 3{\cdot} 5{\cdot} 7,\ 2^2{\cdot} 3^2{\cdot} 7,\ 2^2{\cdot} 3^2{\cdot} 5,\ 2^3{\cdot} 3{\cdot} 7,\ 2^3{\cdot} 3{\cdot} 5,\ 2^3{\cdot} 3^2\}$ \\[1pt] \hline

 $J_1$ &
$\{1,\ 3{\cdot} 5{\cdot} 11{\cdot} 19,\ 2{\cdot} 3{\cdot} 5{\cdot} 7{\cdot} 11,\ 2^3{\cdot} 3{\cdot} 5{\cdot} 19,\ 7{\cdot} 11{\cdot} 19,\ 2^3{\cdot} 3{\cdot} 5{\cdot} 11,\ 2^3{\cdot} 3{\cdot} 5{\cdot} 7 \}$\\[1pt]  \hline

  $M_{11}$ &
$\{ 1,\ 2^3{\cdot} 3^2{\cdot} 11,\ 2^4{\cdot} 3^2{\cdot} 5,\ 3^2{\cdot} 5{\cdot} 11,\ 2^2{\cdot} 3^2{\cdot} 5,\ 2^4{\cdot} 11,\ 2^4{\cdot} 3^2 \}$ \\ [1pt]\hline

 $\PSL(3,3)$ &
$\{1,\ 2^2{\cdot} 3^2{\cdot} 13,\ 2^4{\cdot} 3^3,\ 3^3{\cdot} 13,\ 2^3{\cdot} 3^3,\ 2^4{\cdot} 13,\ 2^4{\cdot} 3^2\}$ \\ [1pt] \hline

  $\PSL(3,q)$ &
$\{ 1, \  (q^2+q+1)(q^2-1)(q-1), \ q^2(q^2+q+1)(q-1)^2, $ \\ [1pt]
$4 < q \not\equiv 1\ (\bmod\ 3)$ & $ \ q^3(q^2+q+1),\ q^2(q^2-1)(q-1),q^3(q^2-1),$ \\ [1pt]
&$ q^3(q^2-1)(q-1), \ q^3(q-1)^2 \}$
\\ [1pt] \hline

$\PSL(3,q)$ &    $ \{ 1, \  \frac{1}{3}(q^2+q+1)(q+1)(q-1)^2, \ \frac{1}{3} q^2 (q^2+q+1)(q-1)^2, $  \\ [1pt]
$4 < q \equiv 1 \ (\bmod\ 3)$ & $  \frac{1}{3} q^3(q^2+q+1), \frac{1}{3} q^2(q+1)(q-1)^2, \ \frac{1}{3}q^3(q-1)(q+1), $  \\ [1pt]
& $\frac{1}{3} q^3 (q+1)(q-1)^2, \ \frac{1}{3} q^3(q-1)^2, \ q^3(q-1)^2 \}$ \\ [1pt] \hline

 $\PSU(3,q)$ &  $\{ 1, \ (q^2-q+1)(q+1)^2(q-1), \ q^3(q^2-q+1),$ \\ [1pt]
 $4 < q \not\equiv -1\ (\bmod\ 3)$  & $\ q^2(q^2-q+1)(q+1)^2, \ q^3(q+1)^2(q-1), \ q^3(q+1)^2,$\\ [1pt]
 &$  q^2(q+1)^2(q-1), \ q^3(q-1)(q+1) \}$\\ [1pt] \hline

 $\PSU(3,q)$ &$ \{ 1, \ \frac{1}{3}(q^2-q+1)(q+1)^2(q-1),\ \frac{1}{3}q^3(q^2-q+1), $  \\ [1pt]
  $4 < q \equiv -1 \ (\bmod\ 3)$ & $ \frac{1}{3}q^2(q^2-q+1)(q+1)^2,\ \frac{1}{3}q^3(q+1)^2(q-1),\ \frac{1}{3}q^3(q+1)^2,$  \\ [1pt]
   & $\frac{1}{3}q^2(q+1)^2(q-1), \ \frac{1}{3}q^3(q-1)(q+1), \ q^3(q+1)^2 \}$ \\ [1pt] \hline
 \end{tabular}
\end{table}

At the end of this section, we introduce two results about the maximal subgroups of $\SL(3,q)$ and $\SU(3,q)$ which are from page 378 and 379 in \cite{bhr}. Here the notation of maximal subgroups in two tables are from \cite{bhr}.

\begin{lem}\label{maxSL3q}
Let $S$ be a maximal subgroup  of $\mathrm{SL}(3,q)$ where $q$ is a power of prime $p$. Then $S$ is isomorphic
to one of the groups in the following table, where $d=|Z(\mathrm{SL}(3,q))|=\gcd(3,q-1)$.

$$
\begin{array}[Table 1]{| c | c |} \hline
    \mbox{Maximal subgroup structure} & \mbox{Conditions} \\ \hline
    E_q^2 : \GL(2,q) & \\ \hline
     (q-1)^2 : S_3 & q \geq 5 \\ \hline
    (q^2 + q + 1) : 3 & q \neq 4 \\ \hline
    \SL(3,q_0) . \left(\frac{q-1}{q_0 - 1}, 3\right) & q = q_0^r, \mbox{ $r$ a prime} \\ \hline
    3_+^{1 + 2} : Q_8 . \frac{(q-1, 9)}{3} & p = q\equiv 1\ (\bmod\ 3) \\ \hline
    d\times \mbox{SO}(3,q) & q \mbox{ odd} \\ \hline
    (q_0 - 1, 3)\times SU(3,q_0) & q = q_0^2 \\ \hline
    d\times \PSL_2(7)  & q=p\equiv 1, 2, 4\ (\bmod\ 7),\  q \neq 2 \\ \hline
    3^{.}A_6  & q=p\equiv 1,4\ (\bmod\ 15);\  q=p^2,\ p\equiv 2,3\ (\bmod\ 5),\ p\neq 3 \\ \hline
\end{array}$$

\end{lem}

\vskip 3mm

\begin{lem}\label{maxSU3q}
Let $S$ be a maximal subgroup  of $\mathrm{SU}(3,q)$ where $q$ is a power of prime $p$. Then $S$ is isomorphic to one of the groups in the following table, where $d=|Z(\mathrm{SU}(3,q))|=\gcd(3,q+1)$.

$$
\begin{array}{| c | c |} \hline
    \mbox{Maximal subgroup structure} & \mbox{Conditions} \\ \hline
    E_q^{1+2} : (q^2-1) & \\ \hline
    GU_2(q) & \\ \hline
     (q+1)^2 : S_3 & q \geq 5 \\ \hline
    (q^2 - q + 1) : 3 & q \neq 3,5 \\ \hline
    \mbox{SU}(3,q_0) . \left(\frac{q+1}{q_0 + 1}, 3\right) & q = q_0^r, \mbox{ $r$ odd prime} \\ \hline
    3_+^{1 + 2} : Q_8 . \frac{(q+1, 9)}{3} & p = q\equiv -1\ (\bmod\ 3), \, q \,\geq 11 \\ \hline
    d\times \mbox{SO}(3,q) & q \mbox{ odd},\ q\geq 7 \\ \hline
    d\times \PSL_2(7)  & q=p\equiv 3, 5, 6\ (\bmod\ 7),\  q \neq 5 \\ \hline
    3^{.}A_6  & q=p\equiv 11,14\ (\bmod\ 15) \\ \hline
     3^{.}{A_6}^.{2_3}  & q=5 \\ \hline
      3^{.}A_7  & q=5 \\ \hline
\end{array}
$$
\end{lem}

\vskip 3mm

\begin{remark}
It follows from \cite[Lemma 3.6]{Aliz}, that if $\cod(G) = \cod(H)$
where $H$ is a non-abelian simple group, then $G$ is perfect.
\end{remark}

\vskip 3mm

\section{Main Result for PSL(3,$q$)}

For an integer $n$ and a prime $p$, we denote by $n_p$ the $p$-part of $n$.

\begin{lem} \label{pslquotient1}
Let $G$ be a finite group with $\cod(G)=\cod(\PSL(3,q))$, where $4 < q \not\equiv 1\ (\bmod\ 3)$. If $N$ is a maximal normal subgroup of $G$, then $G/N \cong \PSL(3,q)$.
\end{lem}

\begin{proof} Let $N$ be a maximal normal subgroup of $G$. Since $\cod(G)=\cod(\PSL(3,q))$, we see that $G$ is perfect. Then $G/N$ is a non-abelian simple group. Since $\cod(G/N)\subseteq \cod(G)$, we see that $|\cod(G/N)|$ is either $4, 5, 6, 7$, or $8$.\\

Suppose $|\cod(G/N)| = 4$. Then $G/N\cong \PSL(2,k)$ where $k=2^f\geq 4$. Then $\cod(G/N)=\{1, k(k-1), k(k+1), k^2-1\}$.
Suppose $q$ is even.
Looking at the 2-part of $k(k-1)$ and $k(k+1)$ we see that either $k=q^2$ or $k=q^3$. However, neither $q^4-1$ nor $q^6-1$ is in $\cod(G)$.
Suppose  $q$ is a power of odd prime $r$. Then $q^3(q^2+q+1) = k^2-1$ since these are the only nontrivial odd codegrees in each set.  If $r \mid k-1$, then $q^3$ divides $k-1$ for $(k-1,k+1)=1$ and $k+1$ divides $q^2+q+1$. We obtain a contradiction with $k + 1 < k - 1$. If $r \mid k+1$, then  $q^3$ divides $k+1$   and $k-1$ divides $q^2+q+1$. We obtain a contradiction with $q^3-(q^2+q+1)=(q-2)(q^2+q+1)+1>2$.\\

 Suppose $|\cod(G/N)|=5$. Then $G/N\cong \PSL(2,k)$ where $k$ is an odd prime power and  $\cod(G/N)=\left\{ 1, \frac{k(k-1)}{2}, \frac{k(k+1)}{2}, \frac{k^2-1}{2}, k(k-\epsilon(k))\right\}$ where $\epsilon(k)=(-1)^{(k-1)/2}$.  Then $k(k-\epsilon(k))/2$,   $k(k-\epsilon(k))\in \cod(G)$, a contradiction for we can't find a codegree is the half of another codegree in $\cod(G)$.\\

 Suppose $|\cod(G/N)| =6$. Then $\cod(G/N) = \cod\left({}^2B_2\left(2^{2f+1}\right)\right)$ or $\cod(\PSL(3,4))$.

 Suppose $\cod(G/N) = \cod\left({}^2B_2(s)\right)$ with $s=2^{2f+1}$ and $r=2^f$. If $q$ is even, then $s^2 = q^3$, as $s^2$ is the largest $2$-part of three codegrees in $\cod\left({}^2B_2(s)\right)$ and $2^{3f+2}=q^2$, as there are only two nontrivial $2$-parts of the codegrees in $\cod(\PSL(3,q))$ and  $\cod\left({}^2B_2(s)\right)$.
 It can be checked that $q=2^f$. Then $s^2 = q^3=2^{3f}$, a contradiction.
 If $q$ is odd, there are three codegrees in $\cod(\PSL(3,q))$ with the same $2$-part, i.e. the $2$-part of $(q+1)(q-1)^2$. Then $s^2(s+2r+1)=q^3(q+1)(q-1)^2$,
  $s^2(s-1)=(q^2+q+1)(q+1)(q-1)^2$ and  $s^2(s-2r+1)=q^2(q+1)(q-1)^2$. We have $s+2r+1=q(s-2r+1)$ which is a contradiction.\\

 Suppose $\cod(G/N) = \cod(\PSL(3,4)) =\{ 1,\  2^4\cdot3^2\cdot7,\ 2^6\cdot3^2,\ 2^6\cdot5,\ 2^6\cdot7,\ 3^2\cdot5\cdot7 \}$. Note that $3^2\cdot 5\cdot 7=315$ is the only nontrivial odd codegree.
 If $q$ is even, $(q^2+q+1)(q+1)(q-1)^2=315$ is a contradiction for $q-1>3$.  If $q$ is odd, $q^3(q^2+q+1)=315$ is a contradiction since $315$ is not divisible by a cube of prime.\\

 Suppose $|\cod(G/N)| = 7$. Suppose $\cod(G/N)= \cod(\PSL(3,3))$. If $q$ is even, then $(q^2+q+1)(q+1)(q-1)^2 = 3^3\cdot 13$, as they are the only nontrivial odd codegrees of each set. From this we obtain that $q \not\in \Z$, a contradiction. If $q$ is odd, then $q^3(q^2+q+1) = 3^3\cdot 13$ by the same reason. Then $q=3$, a contradiction.

 Suppose $\cod(G/N)= \cod(A_7)=\{1,\ 2^2\cdot 3\cdot 5\cdot 7,\ 2^2\cdot 3^2\cdot 7,\ 2^2\cdot 3^2\cdot 5,\ 2^3\cdot 3\cdot 7,\ 2^3\cdot 3\cdot 5,\ 2^3\cdot 3^2\}$. Since every nontrivial codegree in  $\cod(A_7)$ is even, then  $\cod(A_7)$ is the subset of $\cod(G)$ after deleting the unique nontrivial odd codegree. If $q$ is even, then by comparing the 2-parts of the codegrees of $A_7$ and $\cod(G)$, we have that $q^3 = 2^3$ which implies that $q=2$, a contradiction. If $q$ is odd.
 Then $2^3$ is the $2$-part of $(q+1)(q-1)^2$ which means that the $2$-part of
 $q-1$ and $q+1$ will be $2$. We obtain a contradiction.

Suppose $\cod(G/N) = \cod(M_{11})$. If $q$ is even, then $q^3 = 2^4$
for there are three codegrees in $ \cod(M_{11})$ with 2-part $2^4$. This is a contradiction.
 If $q$ is odd, we have $3^2\cdot 5\cdot 11 = q^3(q^2+q+1)$ since they are the only nontrivial odd codegrees. A contradiction since $3^2\cdot 5\cdot 11$ is not divisible by a cube of a prime.

 Suppose $\cod(G/N) = \cod(J_1)$. We need only note that $\cod(J_1)$ has two nontrivial odd codegrees while $\cod(G)$ will only have one nontrivial odd codegree.\\

 So, $|\cod(G/N)|=8$. Suppose $\cod(G/N) = \cod(\PSU(3,f))$ with $4 < f \not \equiv -1\ (\bmod\ 3)$.

 Since $f^3(f-1)(f+1)^2$ is the unique codegree  which is divided by another nontrivial  codegree in $\cod(\PSU(3,f))$, then $q^3(q-1)^2(q+1)=f^3(f-1)(f+1)^2$. Furthermore $q^2(q-1)^2(q+1)=f^2(f-1)(f+1)^2$, $f^3(f+1)^2$ or $f^3(f+1)(f-1)$.
 It can be easily checked that $q=f, f-1$ or $f+1$, a  contradiction with $q^3(q-1)^2(q+1)=f^3(f-1)(f+1)^2$.

 Suppose $\cod(G/N) = \cod(G_2(2)')$.  If $q$ is even, then $(q^2+q+1)(q+1)(q-1)^2 = 3^3\cdot 7$. This has no integer solution. If $q$ is odd, then $q^3(q^2+q+1)=3^3\cdot 7$. This also has no integer solution.

 Thus, $\cod(G/N) = \cod(\PSL(3,f))$ for some $4 < f \not \equiv 1\ (\bmod\ 3)$. Comparing the smallest codegrees, we see that $q=f$. Thus, $G/N\cong \PSL(3,q)$.
\end{proof}

\begin{lem} \label{pslquotient2}
Let $G$ be a finite group with $\cod(G)=\cod(\PSL(3,q))$, where $4 < q \equiv 1\ (\bmod\ 3)$. If $N$ is a maximal normal subgroup of $G$, then $G/N \cong \PSL(3,q)$.
\end{lem}

\begin{proof}
Let $N$ be a maximal normal subgroup of $G$. Since $\cod(G)=\cod(\PSL(3,q))$, we see that $G$ is perfect. Then $G/N$ is a non-abelian simple group. Since $\cod(G/N)\subseteq \cod(G)$, we see that $|\cod(G/N)|$ is either $4, 5, 6, 7, 8$, or $9$.\\

Suppose $|\cod(G/N)| = 4$. Then $G/N\cong \PSL(2,k)$ where $k=2^f\geq 4$. Then $\cod(G/N)=\{1, k(k-1), k(k+1), k^2-1\}$.

 Suppose $q$ is even.
Looking at the 2-part of $k(k-1)$ and $k(k+1)$ we see that either $k=q^2$ or $k=q^3$. However, neither $q^4-1$ nor $q^6-1$ is in $\cod(G)$, a contradiction.

 Suppose $q$ is a power of an odd prime $r$. Then $\frac{1}{3}q^3(q^2+q+1) = k^2-1$, as they are the only nontrivial odd codegrees in each set. If $r\mid k-1$, then $q^3 \mid k-1$ since $(k+1,k-1) = 1$ and $k+1 \mid \frac{1}{3} (q^2 + q + 1)$. We have that $k+1<k-1$,  a contradiction.  If $r \mid k+ 1$, then $q^3 \mid k + 1$ and $k - 1 \mid \frac{1}{3} (q^2 + q + 1)$. We have a contradiction since $q^3 - \frac{1}{3} ( q^2 + q + 1 )$ has a minimum value when $q = 7$ which gives $324 > 2$.\\

Suppose $|\cod(G/N)|=5$. Then $G/N\cong \PSL(2,k)$ where $k$ is an odd prime power and  $\cod(G/N)=\left\{ 1, \frac{k(k-1)}{2}, \frac{k(k+1)}{2}, \frac{k^2-1}{2}, k(k-\epsilon(k))\right\}$ where $\epsilon(k)=(-1)^{(k-1)/2}$.  Then $k(k-\epsilon(k))/2$,  $k(k-\epsilon(k))\in \cod(G)$, a contradiction for we can't find a codegree is the half of another codegree in $\cod(G)$.\\

Suppose $|\cod(G/N)| =6$. Then $\cod(G/N) = \cod\left({}^2B_2\left(2^{2f+1}\right)\right)$ or $\cod(\PSL(3,4))$.

Suppose $\cod(G/N) = \cod\left({}^2B_2(s)\right)$ with $s=2^{2f+1}$ and $r=2^f$. If $q$ is even, then $s^2 = q^3$, as $s^2$ is the largest 2-part of three codegrees in $\cod\left({}^2B_2(s)\right)$  and $2^{3f+2}=q^2$, as there are only two nontrivial 2-parts of codegrees in $\cod(\PSL(3,q))$ and  $\cod\left({}^2B_2(s)\right)$. It can be checked that $q=2^f$. Then $s^2 = q^3=2^{3f}$, a contradiction.

 If $q$ is odd, then there are three codegrees in $\cod(\PSL(3,q))$ with the same 2-part, i.e. 2-part of $(q^2-1)(q+1)$. Then, $s^2 (s+2r+1) = \frac{1}{3} q^3 (q+1)(q-1)^2$, $s^2(s-1) = \frac{1}{3} (q^2 + q + 1)(q^2 - 1)(q-1)$ and $s^2 (s-2r + 1) = \frac{1}{3} q^2 (q^2 - 1)(q-1)$. Then $s+ 2r + 1 = q (s - 2r + 1)$. This is a contradiction.\\

  Suppose $\cod(G/N) = \cod(\PSL(3,4))$. If $q$ is even, then $\frac{1}{3}(q^2+q+1)(q+1)(q-1)^2 = 3^2 \cdot 5 \cdot 7$, as they are the only nontrivial odd codegrees in each set. Then $q =4$, a contradiction.
  If $q$ is odd, then $\frac{1}{3}q^3(q^2+q+1) = 3^2 \cdot 5 \cdot 7$, as they are the only nontrivial odd codegrees in each set. This has no integer solution, a contradiction.\\

Suppose $|\cod(G/N)| = 7$. Then $\cod(G/N)= \cod(\PSL(3,3)), \cod(A_7), \cod(J_1),$ or $\cod(M_{11})$.

Suppose $\cod(G/N) = \cod(\PSL(3,3))$. Again by setting an equality between the only nontrivial odd codegrees of each set, we obtain no integer solution, for both $q$ even and $q$ odd. A contradiction.

Suppose $\cod(G/N) = \cod(A_7)$. If $q$ is even, we compare the largest $2$-parts of each codegree set and obtain that $q^3 = 2^3$, and therefore $q=2$, a contradiction.
If $q$ is odd,  $2^3$ is the 2-part of $(q^2 -1)(q-1)$. Thus, the 2-parts of $q-1$ and $q+1$ equal to $2$. This is a contradiction.

Suppose $\cod(G/N) = \cod(J_1)$. We note that $\cod(J_1)$ has two nontrivial odd codegrees while $\cod(G)$ has only one, a contradiction.

Suppose $\cod(G/N) = \cod(M_{11})$. If $q$ is even, then $q^3 = 2^4$ since there are three codegrees in $\cod(M_{11})$ with $2$-part $2^4$. This is a contradiction.
If $q$ is odd, then $3^2 \cdot 5\cdot 11 = \frac{1}{3} q^3 (q^2+q+1)$ since they are the only nontrivial odd codegrees. It can be checked that there is no integer solution, a contradiction.\\

Suppose $|\cod(G/N)| = 8$. Then $\cod(G/N) = \cod(\PSL(3,f))$ where $4 < f \not\equiv 1\ (\bmod\ 3)$ or $\cod(G/N) = \cod(\PSU(3,f))$ where $4 < f \not\equiv -1 \ (\bmod\ 3)$ or $\cod(G/N) = \cod(G_2(2)')$.

Suppose $\cod(G/N) = \cod(\PSL(3,f))$ where $4 < f \not\equiv 1 \ (\bmod\ 3)$. Then
$f^3(f+1)(f-1)^2=\frac{1}{3}q^3(q+1)(q-1)^2$ for there is only one nontrivial codegree which is divided by other three codegrees in $\cod(G)$ and $ \cod(\PSL(3,f))$ where $4 < f \not\equiv 1\ (\bmod\ 3)$.
Furthermore $f^3 (f-1)^2=\frac{1}{3} q^3 (q-1)^2$ as
$f^3 (f-1)^2$ is the smallest codegree among three codegrees which divide same codegree in $\cod(\PSL(3,f)$ and $\frac{1}{3} q^3 (q-1)^2$ satisfies similar condition in   $\cod(G)$. Hence $f+1=q+1$, a contradiction.

Suppose $\cod(G/N) = \cod(\PSU(3, f))$ where $4 < f \not\equiv -1\ (\bmod\ 3)$. Similar with the analysis in above paragraph, we have  $f^3(f+1)^2(f-1)=\frac{1}{3}q^3(q+1)(q-1)^2$ and $f^3(f+1)(f-1)=\frac{1}{3}q^3(q-1)^2$. Hence $f+1=q+1$, a contradiction.

 Suppose $\cod(G/N) = \cod(G_2(2)')$. If $q$ is even, then $\frac{1}{3}(q^2+q+1)(q+1)(q-1)^2 = 3^3\cdot 7$. This has no integer solution. If $q$ is odd, then $\frac{1}{3}q^3(q^2+q+1)=3^3\cdot 7$. This also has no integer solution.\\

Thus, $|\cod(G/N)| = 9$. Suppose $\cod(G/N) = \cod(\PSU(3,f))$ where $4 < f \equiv -1\ (\bmod\ 3)$.
We have  $\frac{1}{3}f^3(f+1)^2(f-1)=\frac{1}{3}q^3(q+1)(q-1)^2$ and $\frac{1}{3}f^3(f+1)(f-1)=\frac{1}{3}q^3(q-1)^2$. Hence $f+1=q+1$, a contradiction.\\

If $\cod(G/N) = \cod(\PSL(3,f))$ where $4 < f \equiv 1\ (\bmod\ 3)$, then $q = f$ by comparing the smallest codegrees. Thus, $G/N\cong \PSL(3,q)$.
\end{proof}

\begin{thm} \label{firstmain} Let $G$ be a group such that  $\cod(G)=\cod(\PSL(3,q))$. Then $G \cong \PSL(3,q)$.
\end{thm}

\begin{proof}

    Let $G$ be a group with $\cod(G) = \cod(\PSL(3,q))$. Let $N$ be a maximal normal subgroup of $G$. Then, $G/N \cong \PSL(3,q)$ by Lemmas \ref{pslquotient1} and \ref{pslquotient2}. Assume to the contrary that $G$ is a minimal counterexample.
    By the choice of $G$,  $N$ is a minimal normal subgroup of $G$. Otherwise there exists a nontrivial normal subgroup  $L$ of $G$ such that $L$ is included in $N$. Then  $\cod(G/L)=\cod(G)$ for  $\cod(G)=\cod(G/N)\subseteq \cod(G/L) \subseteq\cod(G)$  and $G/L\cong \PSL(3,q)$ for $G$ is a minimal counterexample, a contradiction.\\

    {\bf Step 1:} $N$ is the unique minimal normal subgroup of $G$.

Otherwise we assume $M$ is another proper nontrivial  normal subgroup of $G$. If $N$ is included in $M$, then $M=N$ or $M=G$ for $G/N$ is simple, a contradiction. Then $N\cap M=1$ and $G=N\times M$. Since  $M$ is also a maximal normal subgroup of $G$, we have $N\cong M\cong \PSL(3,q)$. Choose $\psi_1\in \irr(N)$ and $\psi_2\in \irr(M)$ such that $\cod(\psi_1)=\cod(\psi_2)=q^3 (q-1)^2$ . Set $\chi=\psi_1\cdot\psi_2\in \irr(G)$. Then $\cod(\chi)=(q^3 (q-1)^2)^2\notin \cod(G)$, a contradiction.\\

Set $\irr(G|N)=\{\chi\in \irr(G)|\,N$ is not contained in the kernel of
$\chi \}$.\\

{\bf Step 2:} $\chi$ is faithful for each $\chi \in \irr(G|N)$.

Since $N$ is not contained in the kernel of $\chi$ for each $\chi \in \irr(G|N)$, we have that the kernel of $\chi$ is trivial by Step 1.\\

{\bf Step 3:} $N$ is  elementary abelian.

Assume to the contrary that $N$ is not abelian. Since $N$ is non-abelian, $N=S^n$ where $S$ is a non-abelian simple group and $n\in \N$.
By Theorem 2, 3, 4 and Lemma 5 in \cite{bcl},  we see that there exists a non-principal character $\psi\in \irr(N)$ that extends to some $\chi\in \irr(G)$. Then $\ker(\chi)=1$ by Step 2 and $\cod(\chi)=|G|/\chi(1)=|G/N|\cdot|N|/\psi(1)$. This is a contradiction since $|G/N|$ is divisible by $\cod(\chi)$. \\

{\bf Step 4:} It is enough to assume that $C_G(N)=N$.

We first note that ${\bf{C}}_G(N) \unlhd G$. Since $N$ is abelian by Step 3, there are two cases: either ${\bf{C}}_G(N) = G$ or ${\bf{C}}_G(N) = N$. If ${\bf{C}}_G(N) = N$, we are done.

Suppose ${\bf{C}}_G(N) = G$. Therefore $N$ must be in the center of $G$.  Since $G$ is perfect, we must have ${\bf Z}(G)=N$ and $N$ is isomorphic to a subgroup of the Schur multiplier of $G/N$ \cite[Corollary 11.20]{Isaacs}. By \cite{gls}, the Schur multiplier of $\PSL(3,q)$ is cyclic of order $\gcd(3, q-1)$. Then $q\equiv 1\ (\bmod\ 3)$ and $G\cong \SL(3,q)$ and  we can find a new codegree $q^3(q+1)(q-1)^2$ in \cite{sf}, a contradiction.
Thus, $C_G(N)=N$.\\

{\bf Step 5:} Let $\lambda$ be a non-principal character in $\irr(N)$ and $\theta \in \irr(I_G(\lambda)|\lambda)$. We show that
$\frac{|I_G(\lambda)|}{\theta(1)} \in \cod(G)$. Also, $\theta(1)$ divides $|I_G(\lambda)/N|$ and $|N|$ divides $|G/N|$.

Let $\lambda$ be a non-principal character in $\irr(N)$. Given $\theta\in \irr(I_G(\lambda)|\lambda)$.
Note that $\chi=\theta^G\in \irr(G)$ and $\chi(1)=|G:I_G(\lambda)|\cdot\theta(1)$ by Clifford theory (see \cite[chapter 6]{Isaacs}). Then $\ker(\chi)=1$ by Step 2 and $\cod(\chi)=\frac{|I_G(\lambda)|}{\theta(1)}$. Especially, we have that $\theta(1)$ divides  $|I_G(\lambda)/N|$, and then $|N|$ divides
$\frac{|I_G(\lambda)|}{\theta(1)}$.
Since $\cod(G)=\cod(G/N)$ and  $|G/N|$ is divisible by every element in $\cod(G/N)$, we have that $|N|\mid |G/N|$.

Next we show $I_G(\lambda)<G$. Otherwise we may assume $I_G(\lambda)=G$. Then $\mathrm{ker}(\lambda)\unlhd G$. Furthermore $\mathrm{ker}(\lambda)=1$ by Step 1 and $N$ is a cyclic subgroup with prime order by Step 3. Therefore $G/N$ is abelian for $G/N=N_G(N)/C_G(N) \leq\Aut(N)$ by the  Normalizer-Centralizer Theorem, a contradiction.\\

{\bf Step 6:} Final contradiction.

Given a non-principle character $\lambda\in \irr(N)$. Let $T:=I_G(\lambda)$. By Step 5, we have that $\frac{|T|}{\theta(1)}\in \cod(G)$ for all $\theta\in \irr(T|\lambda)$.

Since $N$ is abelian by Step 1, $|\irr(N)|=|N|$. Therefore, $|N|=|\irr(N)| > |G:T|$ since $|G:T|$ is the number of conjugates of $\lambda$ in $G$ which are all contained in $\irr(N)$. Note that $\gcd(q^2 + q + 1, q + 1)= 1$ and $\gcd(q^2 + q + 1, q - 1)= 1$ and $\gcd(q + 1, q - 1) = 1$ or $2$. Thus, $q^3$ is the largest power of a prime that divides the order of $\PSL(3,q)$.
Then, $|N| \leq q^3$.
Let $Z=Z(\SL(3,q))$ and $K$ be any maximal subgroup of $\SL(3,q)$ such that $T/N$ is isomorphic to a subgroup of $K/Z$. Then $K/Z$ is isomorphic to  a maximal subgroup of $G/N$. If $K$ is not of the type $E_q^2 : \GL(2,q)$ in Lemma \ref{maxSL3q}, then $|G:T|>q^3$, a contradiction. Therefore $K$ must be of the type $E_q^2 : \GL(2,q)$ and $|G:T|\geq \frac{|G/N|}{|K/Z|}=q^2+q+1$. Assume $|N|$ is a power of prime $p$, then $|N|=2(q-1)^2$ where $q-1$ is 2-power or   $q^2\mid |N|\mid q^3$ where $q$ is a power of $p$. If $|N|=2(q-1)^2$ where $q-1$ is $2$-power or $|N|=q^3$, we have that $\frac{|T/N|}{\theta(1)}$ is coprime with $p$ for $|N|$ is the largest $p$-part of codegree. Let $a$ be the $p$-part of $|T/N|$, then $a$ is also the the $p$-part of $\theta(1)$. By Clifford theory of characters we have that $|T/N|$ is the sum of  $\theta(1)^2$ for all $\theta\in \irr(T|\lambda)$ (see chapter 6 of \cite{Isaacs} for more details). Then $a$ is divisible by $a^2$. Therefore $a=1$, i.e. $|T/N|$ is coprime with $q$. So $T/N$ is isomorphic to
a subgroup of  $\GL(2,q)/Z$. We obtain a contradiction for $|G:T|>q^3$. Next we assume $q^2\mid |N|$ and $|N|< q^3$ where $q$ is a power of $p$. It can be checked that the $p$-part of $\frac{|T|}{\theta(1)}$ is at most $q^3$  and $(\frac{|T/N|}{\theta(1)})_p<q$, or equivalently, $\frac{|T/N|_p}{q} \leq \theta(1)_p$ for every
$\theta\in \irr(T|\lambda)$. Note that    $|T/N|$ is the sum
of $\theta(1)^2$ for all $\theta\in \irr(T|\lambda)$ . Thus, $|T/N|_p<q^2$. Therefore $|G:T|>q(q^2+q+1)$, a contradiction.
\end{proof}

\section{Main Result For PSU(3,$q$)}

\begin{lem} \label{psuquotient1}
Let $G$ be a finite group with $\cod(G)=\cod(\PSU(3,q))$, where $4 < q \not\equiv -1\ (\bmod\ 3)$. If $N$ is a maximal normal subgroup of $G$, then $G/N \cong \PSU(3,q)$.
\end{lem}

\begin{proof} Let $N$ be a maximal normal subgroup of $G$. Since $\cod(G)=\cod(\PSU(3,q))$, we see that $G$ is perfect. Then $G/N$ is a non-abelian simple group. Since $\cod(G/N)\subseteq \cod(G)$, we see that $|\cod(G/N)|$ is either $4, 5, 6, 7$, or $8$.\\

Suppose $|\cod(G/N)| = 4$. Then $G/N\cong \PSL(2,k)$ where $k=2^f\geq 4$. Then $\cod(G/N)=\{1, k(k-1), k(k+1), k^2-1\}$.
Suppose $q$ is even.
Looking at the 2-part of $k(k-1)$ and $k(k+1)$ we see that either $k=q^2$ or $k=q^3$. However, neither $q^4-1$ nor $q^6-1$ is in $\cod(G)$.
Suppose  $q$ is a power of odd prime $r$. Then $q^3(q^2-q+1) = k^2-1$ since these are the only nontrivial odd codegrees in each set.  If $r \mid k-1$, then $q^3$ divides $k-1$ for $(k-1,k+1)=1$ and $k+1$ divides $q^2+q+1$. We obtain a contradiction. If $r \mid k+1$, then  $q^3$ divides $k+1$   and $k-1$ divides $q^2+q+1$. We obtain a contradiction with $q^3-(q^2-q+1)=(q-1)(q^2+1)>2$.\\

 Suppose $|\cod(G/N)|=5$. Then $G/N\cong \PSL(2,k)$ where $k$ is an odd prime power and  $\cod(G/N)=\left\{ 1, \frac{k(k-1)}{2}, \frac{k(k+1)}{2}, \frac{k^2-1}{2}, k(k-\epsilon(k))\right\}$ where $\epsilon(k)=(-1)^{(k-1)/2}$.  Then $k(k-\epsilon(k))/2$,   $k(k-\epsilon(k))\in \cod(G)$, a contradiction for we can't find a codegree is the half of another codegree in cod(G).\\

 Suppose $|\cod(G/N)| =6$. Then $\cod(G/N) = \cod({}^2B_2(2^{2f+1}))$ or $\cod(\PSL(3,4))$.

 Suppose $\cod(G/N) = \cod({}^2B_2(s))$ with $s=2^{2f+1}$ and $r=2^f$. If $q$ is even, then $s^2 = q^3$, as $s^2$ is the largest 2-part of three codegrees in $\cod({}^2B_2(s))$ and $2^{3f+2}=q^2$, as there are only two nontrivial 2-parts of the codegrees in $\cod(\PSU(3,q))$ and  $\cod({}^2B_2(s))$. It can be checked that $q=2^f$. Then $s^2 = q^3=2^{3f}$, a contradiction.
 If $q$ is odd, there are three codegrees in $\cod(\PSU(3,q))$ with same  2-part, i.e. 2-part of $(q+1)^2(q-1)$. Then $s^2(s+2r+1)=q^3(q+1)^2(q-1)$,
  $s^2(s-1)=q^2(q+1)^2(q-1)$ and  $s^2(s-2r+1)=(q^2-q+1)(q+1)^2(q-1)$. We have $s+2r+1=q(s-1)$ which is a contradiction.\\

 Suppose $\cod(G/N) = \cod(\PSL(3,4)) =\{ 1,\  2^4\cdot3^2\cdot7,\ 2^6\cdot3^2,\ 2^6\cdot5,\ 2^6\cdot7,\ 3^2\cdot5\cdot7 \}$. Note that $3^2\cdot 5\cdot 7=315$ is the only nontrivial odd codegree.
 If $q$ is even, $(q^2-q+1)(q+1)^2(q-1)=315$ is a contradiction for $q+1>3$.  If $q$ is odd, $q^3(q^2-q+1)=315$ is a contradiction since $315$ is not divisible by a cube of prime.\\

 Suppose $|\cod(G/N)| = 7$. Suppose $\cod(G/N)= \cod(\PSL(3,3))$. If $q$ is even, then $(q^2-q+1)(q+1)^2(q-1) = 3^3\cdot 13$, as they are the only nontrivial odd codegrees of each set. From this we obtain that $q \not\in \Z$, a contradiction. If $q$ is odd, then $q^3(q^2-q+1) = 3^3\cdot 13$ by the same reason. Then $q=3$, a contradiction.

 Suppose $\cod(G/N)= \cod(A_7)=\{1,\ 2^2\cdot 3\cdot 5\cdot 7,\ 2^2\cdot 3^2\cdot 7,\ 2^2\cdot 3^2\cdot 5,\ 2^3\cdot 3\cdot 7,\ 2^3\cdot 3\cdot 5,\ 2^3\cdot 3^2\}$. If $q$ is even, then by comparing the $2$-parts of the codegrees of $A_7$ and $\cod(G)$, we have that $q^3 = 2^3$ which implies that $q=2$, a contradiction. If $q$ is odd.
 Then $2^3$ is the $2$-part of $(q+1)^2(q-1)$ which means that the $2$-parts of
 $q-1$ and $q+1$ will be $2$. We obtain a contradiction.

Suppose $\cod(G/N) = \cod(M_{11})$. If $q$ is even, then $q^3 = 2^4$
for there are three codegrees in $ \cod(M_{11})$ with $2$-part $2^4$. This is a contradiction.
 If $q$ is odd, we have $3^2\cdot 5\cdot 11 = q^3(q^2-q+1)$ since they are the only nontrivial odd codegrees. A contradiction since $3^2\cdot 5\cdot 11$ is not divisible by a cube of a prime.

 Suppose $\cod(G/N) = \cod(J_1)$. We need only note that $\cod(J_1)$ has two nontrivial odd codegrees while $\cod(G)$ will only have one nontrivial odd codegree.\\

 So, $|\cod(G/N)|=8$. Suppose $\cod(G/N) = \cod(\PSL(3,f))$ with $4 < f \not \equiv 1\ (\bmod\ 3)$.

 Since $f^3(f+1)(f-1)^2$ is the unique codegree  which is divided by another nontrivial  codegree in $\cod(\PSL(3,f))$, then $q^3(q+1)^2(q-1)=f^3(f+1)(f-1)^2$. Furthermore $q^2(q+1)^2(q-1)=f^2(f+1)(f-1)^2$, $f^3(f-1)^2$ or $f^3(f+1)(f-1)$.
 It can be easily checked that $q=f, f+1$ or $f-1$ which is a  contradiction with $q^3(q+1)(q-1)^2=f^3(f+1)(f-1)^2$.

 Suppose $\cod(G/N) = \cod(G_2(2)')$.  If $q$ is even, then $(q^2-q+1)(q+1)^2(q-1) = 3^3\cdot 7$. This has no integer solution. If $q$ is odd, then $q^3(q^2-q+1)=3^3\cdot 7$. This also has no integer solution.

 Thus, $\cod(G/N) = \cod(\PSU(3,f))$ for some $4 < f \not \equiv -1\ (\bmod\ 3)$. Comparing the smallest codegrees, we see that $q=f$. Thus, $G/N\cong \PSU(3,q)$.
\end{proof}

\begin{lem} \label{psuquotient2}
Let $G$ be a finite group with   $\cod(G)=\cod(\PSU(3,q))$, where $4 < q \equiv -1\ (\bmod\ 3)$. If $N$ is a maximal normal subgroup of $G$, then $G/N \cong \PSU(3,q)$.
\end{lem}

\begin{proof}
Let $N$ be a maximal normal subgroup of $G$. Since $\cod(G)=\cod(\PSU(3,q))$, we see that $G$ is perfect. Then $G/N$ is a non-abelian simple group. Since $\cod(G/N)\subseteq \cod(G)$, we see that $|\cod(G/N)|$ is either $4, 5, 6, 7, 8$, or $9$.\\

Suppose $|\cod(G/N)| = 4$. Then $G/N\cong \PSL(2,k)$ where $k=2^f\geq 4$. Then $\cod(G/N)=\{1, k(k-1), k(k+1), k^2-1\}$.

 Suppose $q$ is even.
Looking at the 2-part of $k(k-1)$ and $k(k+1)$ we see that either $k=q^2$ or $k=q^3$. However, neither $q^4-1$ nor $q^6-1$ is in $\cod(G)$, a contradiction.

If $q$ is a power of an odd prime $r$, then $\frac{1}{3}q^3(q^2-q+1) = k^2-1$, as they are the only nontrivial odd codegrees in each set. If $r\mid k-1$, then $q^3 \mid k-1$ since $(k+1,k-1) = 1$. Then, $k+1 \mid \frac{1}{3} (q^2 - q + 1)$. We have that $k+1<k-1$,  a contradiction.  If $r \mid k+ 1$, then $q^3 \mid k + 1$ and $k - 1 \mid \frac{1}{3} (q^2 - q + 1)$. We have a contradiction since $q^3 - \frac{1}{3} ( q^2 - q + 1 )$ has a minimum value when $q = 5$ which gives $118 > 2$.\\

Suppose $|\cod(G/N)|=5$. Then $G/N\cong \PSL(2,k)$ where $k$ is an odd prime power and  $\cod(G/N)=\left\{ 1, \frac{k(k-1)}{2}, \frac{k(k+1)}{2}, \frac{k^2-1}{2}, k(k-\epsilon(k))\right\}$ where $\epsilon(k)=(-1)^{(k-1)/2}$.  Then $k(k-\epsilon(k))/2$,  $k(k-\epsilon(k))\in \cod(G)$, a contradiction for we can't find a codegree is the half of another codegree in cod(G).\\

Suppose $|\cod(G/N)| =6$. Then $\cod(G/N) = \cod({}^2B_2(2^{2f+1}))$ or $\cod(\PSL(3,4))$.

Suppose $\cod(G/N) = \cod({}^2B_2(s))$ with $s=2^{2f+1}$ and $r=2^f$. If $q$ is even, then $s^2 = q^3$, as $s^2$ is the largest 2-part of three codegrees in $\cod({}^2B_2(s))$  and $2^{3f+2}=q^2$, as there are only two nontrivial 2-parts of codegrees in $\cod(\PSU(3,q))$ and  $\cod({}^2B_2(s))$. It can be checked that $q=2^f$. Then $s^2 = q^3=2^{3f}$, a contradiction.

 If $q$ is odd, then there are three codegrees in $\cod(\PSU(3,q))$ with the same 2-part, i.e. 2-part of $(q+1)^2(q+1)$. Then, $s^2 (s+2r+1) = \frac{1}{3} q^3 (q+1)^2(q-1)$, $s^2(s-1) = \frac{1}{3} q^2(q+1)^2(q-1)$ and $s^2 (s-2r + 1) = \frac{1}{3} (q^2 - q + 1)(q+1)^2(q-1)$. Then $s+ 2r + 1 = q (s -  1)$. This is a contradiction.\\

  Suppose $\cod(G/N) = \cod(\PSL(3,4))$. If $q$ is even, then $\frac{1}{3}(q^2-q+1)(q+1)^2(q-1) = 3^2 \cdot 5 \cdot 7$, as they are the only nontrivial odd codegrees in each set. A contradiction.
  If $q$ is odd, then $\frac{1}{3}q^3(q^2-q+1) = 3^2 \cdot 5 \cdot 7$, as they are the only nontrivial odd codegrees in each set. This has no integer solution, a contradiction.\\

Suppose $|\cod(G/N)| = 7$. Then $\cod(G/N)= \cod(\PSL(3,3)), \cod(A_7), \cod(J_1),$ or $\cod(M_{11})$.

Suppose $\cod(G/N) = \cod(\PSL(3,3))$. Again by setting an equality between the only nontrivial odd codegrees of each set, we obtain no integer solution, for both $q$ even and $q$ odd. A contradiction.

Suppose $\cod(G/N) = \cod(A_7)$. If $q$ is even, we compare the largest 2-parts of each codegree set and obtain that $q^3 = 2^3$, and therefore $q=2$, a contradiction.
If $q$ is odd,  $2^3$ is the 2-part of $(q+1)^2(q-1)$. Thus, the 2-part of $q-1$ and $q+1$ is $2$. This is a contradiction.

Suppose $\cod(G/N) = \cod(J_1)$. We note that $\cod(J_1)$ has two nontrivial odd codegrees while $\cod(G)$ has only one, a contradiction.

Suppose $\cod(G/N) = \cod(M_{11})$. If $q$ is even, then $q^3 = 2^4$ since there are three codegrees in $\cod(M_{11})$ with $2$-part $2^4$. This is a contradiction.
If $q$ is odd, then $3^2 \cdot 5\cdot 11 = \frac{1}{3} q^3 (q^2-q+1)$ since they are the only nontrivial odd codegrees. It can be checked that there is no integer solution, a contradiction.\\

Suppose $|\cod(G/N)| = 8$. Then $\cod(G/N) = \cod(\PSL(3,f))$ where $4 < f \not\equiv 1\ (\bmod\ 3)$ or $\cod(G/N) = \cod(\PSU(3,f))$ where $4 < f \not\equiv -1\ (\bmod\ 3)$ or $\cod(G/N) = \cod(G_2(2)')$.

Suppose $\cod(G/N) = \cod(\PSL(3,f))$ where $4 < f \not\equiv 1\ (\bmod\ 3)$. Then
$f^3(f+1)(f-1)^2=\frac{1}{3}q^3(q+1)^2(q-1)$ for there is only one nontrivial codegree which is divided by other three codegrees in $\cod(G)$ and $ \cod(\PSL(3,f))$ where $4 < f \not\equiv 1\ (\bmod\ 3)$.
Furthermore $f^3 (f-1)^2=\frac{1}{3} q^3 (q+1)(q-1)$ as
$f^3 (f-1)^2$ is the smallest codegree among three codegrees which divide same codegree in $\cod(\PSL(3,f))$ and $\frac{1}{3} q^3 (q-1)^2$ satisfies similar condition in   $\cod(G)$. Hence $f+1=q+1$, a contradiction.

Suppose $\cod(G/N) = \cod(\PSU(3, f))$ where $4 < f \not\equiv -1\ (\bmod\ 3)$. Similar with the analysis in above paragraph, we have  $f^3(f+1)^2(f-1)=\frac{1}{3}q^3(q+1)^2(q-1)$ and $f^3(f+1)(f-1)=\frac{1}{3}q^3(q+1)(q-1)$. Hence $f+1=q+1$, a contradiction.

 Suppose $\cod(G/N) = \cod(G_2(2)')$. If $q$ is even, then $\frac{1}{3}(q^2-q+1)(q+1)^2(q-1) = 3^3\cdot 7$. This has no integer solution. If $q$ is odd, then $\frac{1}{3}q^3(q^2-q+1)=3^3\cdot 7$. This also has no integer solution.\\

Thus, $|\cod(G/N)| = 9$. Suppose $\cod(G/N) = \cod(\PSL(3,f))$ where $4 < f \equiv 1\ (\bmod\ 3)$.
We have  $\frac{1}{3}f^3(f+1)(f-1)^2=\frac{1}{3}q^3(q+1)^2(q-1)$ and $\frac{1}{3}f^3(f-1)^2=\frac{1}{3}q^3(q+1)(q-1)$. Hence $f+1=q+1$, a contradiction.\\

If $\cod(G/N) = \cod(\PSU(3,f))$ where $4 < f \equiv -1\ (\bmod\ 3)$, then $q = f$ by comparing the smallest codegrees. Thus, $G/N\cong \PSU(3,q)$.
\end{proof}

\begin{thm} Let $G$ be a group such that  $\cod(G)=\cod(\PSU(3,q))$. Then $G \cong \PSU(3,q)$.
\end{thm}

\begin{proof}

    Let $G$ be a group with $\cod(G) = \cod(\PSU(3,q))$. Let $N$ be a maximal normal subgroup of $G$. Then, $G/N \cong \PSU(3,q)$ by Lemmas \ref{psuquotient1} and \ref{psuquotient2}. Assume to the contrary that $G$ is a minimal counterexample.
    By the choice of $G$,  $N$ is a minimal normal subgroup of $G$. Similarly with Theorem \ref{firstmain}, the following six steps can be obtained.\\

    {\bf Step 1:} $N$ is the unique minimal normal subgroup of $G$.\\

{\bf Step 2:} $\chi$ is faithful for each $\chi \in \irr(G|N)$. \\

{\bf Step 3:} $N$ is  elementary abelian. \\

{\bf Step 4:} It is enough to assume that $C_G(N)=N$.

Suppose ${\bf{C}}_G(N) = G$. Therefore $N$ must be in the center of $G$.  Since $G$ is perfect, we must have ${\bf Z}(G)=N$ and $N$ is isomorphic to a subgroup of the Schur multiplier of $G/N$ \cite[Corollary 11.20]{Isaacs}. By \cite{gls}, the Schur multiplier of $\PSU(3,q)$ is cyclic of order $\gcd(3, q+1)$. Then $q\equiv -1\ (\bmod\ 3)$ and $G\cong \SU(3,q)$ and  we can find a new codegree $q^3(q+1)^2(q-1)$ in \cite{sf}, a contradiction.
Thus, $C_G(N)=N$.\\

{\bf Step 5:} Let $\lambda$ be a non-principal character in $\irr(N)$ and $\theta \in \irr(I_G(\lambda)|\lambda)$. We show that
$\frac{|I_G(\lambda)|}{\theta(1)} \in \cod(G)$. Also, $\theta(1)$ divides $|I_G(\lambda)/N|$ and $|N|$ divides $|G/N|$. \\

{\bf Step 6:} Final contradiction.

By Step 3, $N$ is an elementary abelian $r$-subgroup for some prime $r$ and we assume $|N|=r^n$, $n\in \N$. By the Normalizer-Centralizer Theorem, we see that $n > 1$.

Given a non-principle character $\lambda\in \irr(N)$. Let $T:=I_G(\lambda)$. By Step 5, we have that $\frac{|T|}{\theta(1)}\in \cod(G)$ for all $\theta\in \irr(T|\lambda)$.

Since $N$ is abelian by Step 1, $|\irr(N)|=|N|$. Therefore, $|N|=|\irr(N)| > |G:T|$ since $|G:T|$ is the number of conjugates of $\lambda$ in $G$ which are all contained in $\irr(N)$. It can be checked that $q^3$ is the largest power of a prime that divides the order of $\PSU(3,q)$.
Then, $|N| \leq q^3$.
Let $Z=Z(\SU(3,q))$ and $K$ be any maximal subgroup of $\SU(3,q)$ such that $T/N \leq K/Z$.
By Lemma \ref{maxSU3q}, we have that the index of  $K$ in $\SU(3,q)$ is bigger than $q^3$. Therefore $|G:T|>q^3$, a contradiction with $|G:T|<|N|$.
\end{proof}

\section{Acknowledgements}
Liu was supported by NSFC (Grant Nos. 11701421 and 11871011). Yang was supported by grants from the Simons Foundation (\#499532 and \#91809). The collaborative research was also funded by the American Mathematical Society's Ky and Yu-Fen Fan Travel Grant Program.

\bigskip

\noindent \textbf{Data availability Statement:} Data sharing not applicable to this article as no datasets were generated or analysed during the current study.

\bigskip

\noindent \textbf{Competing interests:} The authors declare none.

\end{document}